\documentclass[12pt,a4paper,reqno]{amsart}
\newtheorem{lemma}{Lemma}

\newtheorem{theorem}[lemma]{Theorem}
\newtheorem{corollary}[lemma]{Corollary}
\newcommand{\R}{{\bf R}}
\newcommand{\C}{{\bf C}}
\newcommand{\rme}{{\rm e}}
\newcommand{\rmd}{{\rm d}}
\newcommand{\alp}{\alpha}
\newcommand{\bet}{\beta}
\newcommand{\gam}{\gamma}
\newcommand{\lam}{\lambda}

\newcommand{\f}{\varphi}
\newcommand{\e}{\varepsilon}
\newcommand{\ta}{\tau}
\newcommand{\lap}{{\Delta}}

\newcommand{\Spec}{{\rm Spec}}

\newcommand{\norm}{\Vert}
\renewcommand{\Re}{{\rm Re}}
\renewcommand{\Im}{{\rm Im}}

\newcommand{\Num}{{\rm Num}}

\newcommand{\Schrodinger}{Schr\"odinger }
\newcommand{\gra}{\nabla}
\newcommand{\dpa}{\partial}
\newcommand{\EV}{eigenvalue}
\newcommand{\EVs}{eigenvalues}

\newcommand{\exg}{{\rm e}^{z x}}
\newcommand{\exd}{{\rm e}^{-z x}}
\newcommand{\inte}{\int_{\R}}
\newcommand{\intn}{\int_{\R^N}}
\begin{document}
\title{BOUNDS ON COMPLEX EIGENVALUES AND RESONANCES}
\author{A.A. Abramov, A. Aslanyan and E.B. Davies}
\address{Department of Numerical Methods\\
Computing Centre of RAS\\
Vavilova 40, Moscow 117967, Russia\\
Department of Mathematics\\
King's College London\\
Strand, London WC2R 2LS, UK}
\email{alalabr@ccas.ru \\
aslanyan@mth.kcl.ac.uk \\
e.brian.davies@kcl.ac.uk}
\thanks{We would like to thank H. Siedentop for
several valuable discussions. We also thank Engineering
and Physical Sciences Research Council for support under grant No.
GR/L75443, the Royal Society for funding under an ex-agreement FSU programme,
and the Russian Foundation for Basic Research for support under grant No. 99-01-00331.}
\date{24 November 1999}
%%%%%%%%%%%%%%%%%%%%%
\begin{abstract}
We obtain bounds on the complex eigenvalues of
non-self-adjoint \Schrodinger operators with complex
potentials, and also on the complex resonances of
self-adjoint \Schrodinger operators. Our bounds are
compared with numerical results, and are
seen to provide useful information.
\end{abstract}
%%%%%%%%%%%%%%%%%%%%%%%
\maketitle
\section{Introduction}
\label{1}

We consider self-adjoint and non-self-adjoint \Schrodinger
operators of the form
\[
Hf:=-\lap f+Vf
\]
acting in $L^{2}(\R^{N})$,
where $V$ is a potential vanishing sufficiently rapidly as $|x|\to\infty$.
When $V$ is a complex-valued potential we obtain bounds on
the location of the eigenvalues of $H$, while in the case
when $V$ is real-valued we obtain bounds on its resonances. We
also describe and implement methods for computing eigenvalues and
resonances, and compare them both with each other and with others.

There appears to be rather little literature on either of
these problems. The distribution of the eigenvalues of
self-adjoint \Schrodinger operators has been intensively
studied for several decades, but most of the techniques are
not applicable to non-self-adjoint operators. Similarly most
of the literature concerning resonances is either
computational or considers the semi-classical limit and
resonances close to the real axis. Our results are of quite a
different character.

The principal technique which we use to obtain bounds on \EVs \ and resonances
in sections 2 and 4 involves obtaining
constraints on the numerical range of the operator $H$ and of
others closely related to it. In section 3, however, we
describe a quite different technique for ordinary differential operators
which yields similar bounds but also enables us to prove
that the number of eigenvalues of
such operators is finite for potentials which decay rapidly
enough. Section 5 is devoted to obtaining a first order perturbation formula
for resonances of a \Schrodinger operator. For a general
perturbation theory of resonances see, for instance,
\cite{Ag} and references therein.

The remaining two sections concentrate on the computation of
resonances. We have investigated such problems using a well-known
complex scaling method
for Schr\"odinger operators with dilation
analytic potentials elsewhere \cite{AD}, but here we consider Schr\"odinger
operators in one dimension with rapidly
decaying potentials which need not be analytic. In section 6 we
describe two different methods of computing the resonances of such
operators.
In the final section 7 of the paper we apply our results to
several examples. In particular, for comparison purposes we take
operators for which the location of the resonances has already been determined.
One of the examples of section 7 dates back to the series of papers
\cite{Br,KLM,REB}, another to our earlier paper \cite{AD}. We
find that the two methods of section 6 have similar accuracy, but
these methods are sometimes a lot much useful than the complex
scaling method and sometimes considerably less accurate. In all
cases the key numerical issue is the accurate computation of highly
oscillatory integrals.
%%%%%%%%%%%%%%%%%%%%%%%%%%%%%%
\section{Bounds on eigenvalues using complex scaling }
\label{2}

We initially assume for simplicity
that $V$ is a bounded complex-valued potential
such that $V(x)\to 0$ as $|x|\to\infty$. Since $V$ is then a
relatively compact perturbation of $-\lap$ it follows
routinely from \cite{Kato}
that the spectrum of $H$ consists of
$\R^{+}:=[0,\infty)$  together with a finite or countable number of
eigenvalues of finite multiplicity which can only
accumulate at points of $\R^{+}$.

We also assume that $z\to V(zx)$ is complex analytic for
all $x\in \R^{N}$ and all $z$ in the closed sector
\[
\{ z\in\C: 0\leq 2\arg(z) \leq \alp \}
\]
where $0<\alp\leq\pi$ and $V(zx)\to 0$ as $|x|\to 0$ uniformly
with respect to
$z$ in compact subsets of the sector. We then put
\[
c(\theta):=\norm V(\rme^{i\theta/2}\,\cdot)\norm_{\infty}
\]
for all $\theta$ satisfying $0\leq \theta \leq \alp$. It follows
from the three lines lemma that $c(\theta)$ is a
logarithmically convex function of $\theta$. In the theorem
below and elsewhere we adopt the notation
\[
B(a,r):=\{ z\in \C :|z-a|\leq r\}.
\]

\begin{theorem}
Under the above assumptions on $V$ the eigenvalues $\lam$ of $H$
which do not satisfy $-\alp \leq \arg(\lam)\leq 0$ all lie
in the compact convex set
\[
\bigcap_{0\leq \theta\leq\alp}
\left\{   \rme^{-i\theta}\R^{+}+B(0,c(\theta))  \right\}
\]
and their only possible accumulation point is at $0$.
\label{th1}
\end{theorem}

\begin{proof}
If we put
\[
\hat H_{z}:=-z^{-2}\lap +V(z\,\cdot)
\]
then $\hat H_{z}$ is unitarily equivalent to $H$ for positive
real $z$, so its eigenvalues are independent of $z$. The same
applies by analytic continuation for complex $z$, in a
sense which is familiar in the theory of resonances
\cite{CFKS,HS}. When
\[
H_{\theta}:=-\rme^{-i\theta}\lap +V(\rme^{i\theta/2}\,\cdot)
\]
for $0\leq\theta\leq\alp$
then the essential spectrum  of $H_{\theta}$ depends on
$\theta$, being equal to
$\rme^{-i\theta}\R^{+}$. Since the eigenvalues can only
accumulate in the essential spectrum and the only common
point of the essential spectra of all of these operators is $0$,
this is the only possible limit point of the eigenvalues.

In order to obtain bounds on the position of the
eigenvalues, we need to introduce the numerical range
\[
\Num(\theta)=\{ \langle H_{\theta} f,f\rangle:
\norm f\norm =1 \}\,,
\]
where we assume that all $f$ lie in the $\theta$-independent
domain of $H_{\theta}$. In the current situation $\Num(\theta)$
is a convex set which contains the entire spectrum of
$H_{\theta}$ and is contained in
\[
\rme^{-i\theta}\R^{+}+B(0,c(\theta)).
\]
It follows that if $\lam$ is an eigenvalue of $H$ then
$\lam$ lies in this set
for all $0\leq \theta \leq\alp$.
\end{proof}

\underbar{Note} If the sector of analyticity of $V$ contains
the positive real axis in its interior, then under similar
assumptions one concludes that the only possible
accumulation point of the set of {\it all} eigenvalues
is $0$.

Other bounds of a similar type are available, depending on what
information about the potential one wishes to use.

\begin{theorem}
Under the above assumptions on $V$ the eigenvalues $\lam$ of $H$
which satisfy $\Im(\lam)> 0$ all lie
in the closed convex set
\[
\bigcap_{0\leq \theta\leq\alp}
\left\{ x+iy: x\sin(\theta)+y\cos(\theta)\leq a(\theta)
\right\}\,,
\]
where the convex function $a(\cdot)$ on $[0,\alp]$ is
defined by
\[
a(\theta):=\sup\left\{
\Im(\rme^{i\theta}V(\rme^{i\theta/2}v)):v\in\R^N \right\}\,.
\]
\label{th2}
\end{theorem}

\begin{proof}
The main statement of the theorem depends upon the fact that
\[
\Num(\theta)\subseteq
\left\{ x+iy: x\sin(\theta)+y\cos(\theta)\leq a(\theta)  \right\}\,.
\]
We now explain why $a(\cdot)$ is convex. If we define
\[
g(z,v):=\Im \lbrace z^2V(zv)\rbrace
\]
for $0\leq 2\arg(z)\leq\alp$, then as the real part of an analytic
function $g(\cdot,v)$ is harmonic. Hence
\[
h(z):=\sup \lbrace g(z,v):v\in\R^N \rbrace
\]
is subharmonic and constant on each ray through the origin. This
forces $a(\cdot )$ to be convex.
\end{proof}

It is clear that similar bounds can be obtained for certain
unbounded potentials. In the next theorem we
state a general result of this type which uses quadratic form techniques.
Here we make the same assumptions on the
analyticity of $V$ in the sector as before, but replace the
regularity assumptions as indicated.

\begin{theorem}
Suppose that $V_{\theta}:=V(\rme^{i\theta/2}\,\cdot)$ is
relatively compact with respect to $-\lap$ for all
$0\leq \theta \leq\alp$, and that for all such $\theta$ the
operator inequality
\[
|V_{\theta}|\leq\gam_{\theta}(-\lap) +\bet_{\theta}
\]
holds where $0\leq\gam_{\theta}<1$. Then every
eigenvalue $\lam$ of $H$ which does not satisfy $-\alp\leq\arg(\lam)\leq 0$
lies in the closed convex set
\[
\bigcap_{0\leq \theta \leq\alp}
\left\{  C_{\theta}+B(0,\bet_{\theta})  \right\}\,,
\]
where
\[
C_{\theta}:=\left\{  z:-\theta-\arcsin(\gam_\theta)\leq \arg(z)
\leq  -\theta+\arcsin(\gam_\theta)  \right\}\,.
\]
\label{th3}
\end{theorem}

\begin{proof}
 The only new component of the proof is the
observation that if we put $s:=\langle -\lap f,f\rangle$ then
we obtain
\begin{eqnarray*}
        \Num(H_{\theta}) & \subseteq & \bigcup_{s\in\R^+}
        B(\rme^{-i\theta}s,\sin(\gam_{\theta})s+\bet_{\theta})  \\
         & = & C_{\theta}+B(0,\bet_{\theta}).
\end{eqnarray*}
\end{proof}

%%%%%%%%%%%%%%%%%%%%%%%%%%%%%%%
\section{Resolvent bounds on complex eigenvalues}
\label{3}

This section describes a completely different assumption
under which one can obtain bounds on the complex
eigenvalues, at least in one dimension. We replace the
analyticity assumption on the potential $V$ by the
assumption that $V$ lies in $L^{1}(\R)\cap L^2(\R)$. This is
known to imply that $V$ is a relatively compact perturbation
of $-\lap$ acting in $L^{2}(\R)$, so
the spectrum of $H=-\lap +V$ is
again equal to $\R^{+}$  together with eigenvalues which
can only accumulate on the non-negative real axis. Moreover the
domain of $H$ is equal to $W^{2,2}(\R)$, which is contained in $C_0(\R)$.

\begin{theorem}
Under the above assumptions every eigenvalue $\lam$ of $H$ which
does not lie on the positive real axis satisfies
\[
|\lam|\leq \norm V\norm_{1}^{2}/4.
\]
\label{th4}
\end{theorem}

\begin{proof}
Let $\lam=-z^2$ be an eigenvalue of $H$ where $\Re\,z>0$ and let $f$
be the corresponding eigenfunction, so that $f\in C_0(\R)$. Then
\begin{equation}
(-\lap +z^2)f=-Vf
\label{evp}
\end{equation}
so
\[-f=(-\lap+z^2)^{-1}Vf.
\]
Putting $X:=|V|^{1/2}$, $W:=V/X$ and $g:=Wf\in L^2$, we deduce
that
\begin{equation}
-g=W(-\lap+z^2)^{-1}X g
\label{inte}
\end{equation}
so
\[
-1\in\Spec(W(-\lap+z^2)^{-1}X).
\]
This approach is similar to the Birman--Schwinger principle
for finding \EVs \ of self-adjoint operators (see, for instance, \cite{BS}).
By evaluating the Hilbert-Schmidt norm of the above operator, whose
kernel is
\[
W(x) \frac{\rme^{-z|x-y|}}{2z} X(y)
\]
we deduce that
\[
\int_{\R^2} |W(x)|^2  \frac{\rme^{-2\Re\,z|x-y|}}{4|z|^2}  |X(y)|^2
\rmd x\rmd y \geq 1.
\]
This implies $\norm V\norm_1^2/4|z|^2 \geq 1$.
\end{proof}

It appears more difficult to prove that the number of
eigenvalues of $H$ is finite or to obtain bounds on the
number. The following theorem which treats the operator $H$
in one dimension provides a result of this kind.

\begin{theorem}
Let the potential $V$ satisfy
\begin{equation}
\|V(x) \rme^{\gam x}\|_1  < \infty
\label{dec}
\end{equation}
for all $\gam\in \R$. Then the number of the \EVs \ of $H$ is finite and
all the \EVs \ (including real ones) satisfy
$$
|\lam| \leq \frac{9}{4} \|V\|_1^2.
$$
\label{th5}
\end{theorem}

\begin{proof}
Along with (\ref{evp}) consider the equation
$$
u'' - z^2 u = 0.
$$
Its linearly independent solutions are
$$
u_1(x) = \frac{\sinh(zx)}{z}, \qquad u_2(x) = \cosh(zx).
$$
We shall seek the solution $f(x)$ of (\ref{evp}) in the form
$$
f(x) = C_1(x) u_1(x) + C_2(x) u_2(x),
$$
where the functions $C_1,\,C_2$ satisfy
$$
C'_1(x) = V(x) f(x)\cosh(zx),
\qquad C'_2(x) = -V(x) f(x)\frac{\sinh(zx)}{z}.
$$
We take $C_1(-\infty)=z,\,C_2(-\infty)=1$ which corresponds
to $f(x)\sim \rme^{z x},\,x\to-\infty$. Then
$$
C_1(x) =  z + \int_{-\infty}^{x} V(t) f(t) \cosh(zt)\,\rmd t,
$$
$$
C_2(x) =  1-\int_{-\infty}^{x} V(t) f(t) \frac{\sinh(zt)}{z}\,\rmd
t.
$$
It follows from (\ref{dec}) that
$$
\lim_{x\to+\infty} |C_k(x)| < \infty, \qquad k=1,2.
$$

We define the function
\begin{equation}
\f(z)= C_1(\infty) + z C_2(\infty).
\label{phi}
\end{equation}
As is easily seen,
\begin{equation}
\f(z)= 2z + \inte V(t) \tilde{f}(t)\rmd t,
\label{phi1}
\end{equation}
where
$$
\tilde{f}(x)=f(x) \exd =  1 + \frac{1}{2z}
\int_{-\infty}^{x}
\left(1-\rme^{-2z(x-t)}\right)V(t)\tilde{f}(t) \rmd t .
$$

We observe that if $\lam=-z^2$ is an \EV \ of $H$, then $\f(z)=0,\,\Re\,z \geq 0$.
For $z$ in the right-hand half-plane we have
$$
\left| 1 - \rme^{-2z \tau}\right| \leq 2 \qquad \forall
\tau >0
$$
and, hence
$$
|\tilde{f}(x)| \leq 1 + \frac{1}{|z|}
 \int_{-\infty}^{x} |V(t) \tilde{f}(t)| \rmd t.
$$
From this formula we get
$$
\sup_{ x \in \R} |\tilde{f}(x)| \leq
\frac{|z|}{|z| - \|V\|_1}.
$$
This together with (\ref{phi1}) implies
$$
\f(z) = 2z + \psi(z), \qquad |\psi(z)| \leq
\frac{|z| \|V\|_1}{|z| - \|V\|_1}.
$$
Therefore $\f(z) \not= 0$ for sufficiently large
$|z|$. In particular, if $|z| > \frac{3}{2}\|V\|_1$
then $|\psi(z)| < 2|z|$ and $\f(z)\not= 0$
for all such $z$.

Since $\tilde{f}(x;z)$ is analytic in $z$, the function
$\f$ defined by (\ref{phi}) is entire.
Therefore the number of its zeros in any bounded domain
is finite. The \EVs \ of $H$ can only be found
among the squares of the zeros of $\f(z)$, hence the final result.
\end{proof}

\underbar{Note}
Bounds similar to those in the above two theorems can be obtained for the
corresponding operator on the half-line subject to the boundary condition
$a f(0) + b f'(0) = 0,\,|a|+|b|>0$, by minor adjustments to the proof.

%%%%%%%%%%%%%%%%%%%%%%%%%%%%%%%
\section{Complex scaling and resonance bounds}
\label{4}

We now apply the above results to the analysis of the
complex resonances of a self-adjoint \Schrodinger operator
$H$.
We assume that $V$ is real-valued, and that it can be
analytically continued to the closed sector
$\lbrace z:-\alp\leq \arg(z)\leq 0\rbrace$
with the same bounds on $V_\theta$ as in theorem~\ref{1}. $H$ itself
has no complex eigenvalues but well-known
arguments \cite{CFKS,HS} show that those eigenvalues of $H_{\theta}$
which lie within the sector
\[
\{z:-\theta <\arg (z) \leq 0\}
\]
are independent of $\theta$. New eigenvalues may be
emerge from the essential spectrum $\rme^{-i\theta}\R^+$
as $\theta$ increases. The set
of all such eigenvalues as $\theta$ varies over the
permitted range is called the set of resonances of $H$. It
is known that this definition is equivalent to others which
do not depend upon complex scaling.

We assume that
$V_{\theta}$ is bounded and vanishes at infinity for
all $\theta$ in the permitted range and put
\[
a(\theta):=\sup\left\{
\Im(\rme^{i\theta}V(\rme^{i\theta/2}v)):v\in\R^N \right\}
\]
as before.

\begin{theorem}
Under the above assumptions those resonances of $H$ which
satisfy $-\alp < \arg(z) \leq 0$ all lie within
the closed convex set
\[
S:=\bigcap_{0\leq \theta \leq\alp}
\lbrace  x+iy: x\sin(\theta)+y\cos(\theta)\leq a(\theta)  \rbrace.
\]
\label{th6}
\end{theorem}

\begin{proof}
Given $\theta$, every resonance in the stated sector either has
not yet been uncovered, in which case $-\alp\leq\arg(z)\leq
-\theta$, or it has been uncovered, in which case
\[
x\sin(\theta)+y\cos(\theta)\leq a(\theta).
\]
The proof is completed by a geometrical argument. The set $S$ is
the intersection of half-planes, and hence is closed and convex.
\end{proof}

The boundary of the set $S$ is the envelope of the lines
\[
x\sin(\theta)+y\cos(\theta)= a(\theta).
\]
which is given parametrically by
\begin{eqnarray*}
x&=& a(\theta )\sin(\theta)+ a'(\theta)\cos(\theta)\\
y&=& a(\theta) \cos(\theta)- a'(\theta)\sin(\theta).
\end{eqnarray*}

We use this formula to find the form of the envelope as it emerges
from the real axis.

\begin{theorem}
If $a(\theta)=a_1\theta +a_2\theta^2+O(\theta^3)$ for small $\theta >0$
then $a_1\geq 0$ and $a_2\geq 0$. Assuming $a_2 >0$,
the form of the envelope of $S$ is
\[
y=-\frac{(x-a_1)^2}{4a_2}+o((x-a_1)^2)
\]
for $(x-a_1)$ small and positive.
\label{th7}
\end{theorem}

\begin{proof}
An application of Taylor's theorem shows that
\[
a_1=\max \{ V(x)+x\cdot\gra V(x)/2:x\in\R^N \}
\]
and this implies that
\[
a_1\geq M:=\max\{V(x):x\in\R^N\}\geq 0.
\]
The inequality $a_2 \geq 0$ follows from the fact that the function
$a(\cdot)$ is convex.
The main statement of the theorem depends
upon eliminating $\theta$ from the formulae
\begin{eqnarray}
x&=& a_1+2a_2\theta+O(\theta^2)\nonumber\\
y&=&-a_2\theta^2+O(\theta^3).\label{pf}
\end{eqnarray}

\end{proof}

\underbar{Note} While one expects the hypothesis $a_2>0$ of the above theorem to
hold generically, the example of subsection~\ref{7.2} has $a_2=0$,
and a higher order expansion would be needed to determine the form
of the envelope near the real axis.

Physical intuition suggests that there cannot be resonances $\lam$
with $\Im(\lam)$ very small and $\Re\,\lam>M$, which is confirmed by
the example of subsection~\ref{7.2} (see figure~1), where $\Re\,\lam<M\approx36.79$.
However, theorem~\ref{th7} fails to prove quite as much. In the case in which
$V(x)=x^2\rme^{-x^2/b^2}$ for some $b>0$, a simple computation shows that $a_1/M\sim
1.256$.

A threshold for the resonances of an operator $H$ is defined
to be a real number $\gam$ such
that every resonance $\lam$ of $H$ satisfies
$\Re\,\lam\leq \gam$.

\begin{corollary}
If $\alp \geq \pi/2$ then the resonances of the operator $H$ possess a threshold.
\end{corollary}

\begin{proof}
On choosing $\theta=\pi/2$, we see that
every resonance satisfies $\Re\,\lam\leq a(\pi/2)$.
\end{proof}

\begin{corollary}
Let $r=1$ or $2$ and let
\[
V(x)=\int_{\R^+} \rme^{-s|x|^{r}}\mu(\rmd s)
\]
for all $x\in\R^{N}$, where $\mu$ is a finite signed
measure on $\R^+$. Then
\[
\gam=\int_{\R^+} |\mu|(\rmd s)
\]
is a threshold for the resonances of the operator $H$.
\end{corollary}

%%%%%%%%%%%%%%%%%%%%%%%%%%%%%%%
\section{Perturbation formulae for resonances}
\label{5}

To obtain perturbation formulae let us assume that $V(x)$ is real-valued
and positive and rewrite the integral equation (\ref{inte}) as
$$
- g(x) = \intn W(x) G(x,y;z) W(y) g(y) \rmd y,
$$
where for all $x,\,y$ the Green's function $G(x,y;z)$ is analytic in $z,\,\Re\,z>0$, positive for
$z>0$ and satisfies $G(x,y;\bar{z}) = \bar{G}(x,y;z)$. For certain values of $N$ the kernel of
(\ref{inte}) is known; the argument below does not use the explicit form of $G$.
The above formula defines the operator
$$
A(z) g = \intn W(x) G(x,y;z) W(y) g(y) \rmd y + g.
$$
The resonances therefore are defined as such values of $\lam=-z^2$
that 0 is an \EV \ of $A(z)$. The operator $A(z)-I$ has a symmetric kernel for
all $z,\,\Re\,z \geq 0$, and is self-adjoint
for positive real $z$. However, since we want to consider complex $z$ let us
decompose the operator $A$ and the function $g$ as
$$
A=A_1 + i A_2, \qquad g = g_1 + i g_2,
$$
where the real and imaginary parts $A_{1,2}$ of $A$ are self-adjoint, and $g_{1,2}$ are real-valued.
The (non-linear in $z$) problem
\begin{equation}
A(z) g =0
\label{sa}
\end{equation}
is equivalent to the equation
\begin{equation}
\left( \begin{array}{cc}
-A_2 & A_1 \\
A_1 & A_2
\end{array} \right)
\left( \begin{array}{c}
g_2 \\
g_1
\end{array} \right) = 0.
\label{ma}
\end{equation}
The problem (\ref{ma}) involves a family of {\em self-adjoint} operators so
we can now
use standard perturbation theory techniques valid for such operators.

Suppose we have a resonance $\lam_0=-z_0^2$, i.e. if we put $z=z_0$ in (\ref{sa})
the problem has a solution $g_0$. Consider a perturbed operator $\hat{A}(z,\e)=A(z) + \e B(z) + O(\e^2)$
obtained, for instance, if we replace $V$ by $V_\e=V+\e V_1+ O(\e^2)$ in (\ref{evp})
(examples are to be found in section 7). Let us find the relevant resonance $\lam=-z^2$ assuming that $z =
z_0+O(\e)$. The latter assumption is justified by the above reasoning and by the
fact that the operator function $A(z)$ is analytic in $z$ when $z\not=0$.
It is therefore expanded in powers of $(z-z_0)$: $A(z) = A(z_0) + A'(z_0)(z-z_0) + \ldots$; the same is true for $B(z)$.
Following the conventional approach (cf. \cite{Kato}) we seek a solution of the problem $\hat{A} g = 0$
in the form
$$
(A(z) + \e B(z)) (g_0 + \e g_1) = O(\e^2).
$$
Denote $z-z_0=\e \nu + O(\e^2)$.
Then the above equation can be rewritten as
$$
(A(z_0) + \e A'(z_0)\nu + \e B(z_0)) (g_0 + \e g_1) = O(\e^2)
$$
yielding
\begin{equation}
\nu=-\frac{(B(z_0) g_0,g_0)}{(A'(z_0) g_0,g_0)}
\label{nu}
\end{equation}
(here we use the property $Ag = 0 \ \Rightarrow \ A^*\bar{g}=0$). The operators
involved are easily computed as
$$
A'(z) g = \intn W(x) \frac{\dpa G(x,y;z)}{\dpa z} W(y) g(y) \rmd y,
$$
$$
B(z) g = \intn V_1(x) G(x,y;z) W(y) g(y) \rmd y \ +
$$
$$
\intn W(x) G(x,y;z) V_1(y) g(y) \rmd
y.
$$

Finally, we formulate the main result of this section.
\begin{theorem}
Let $\lam_0$ be a resonance of $H=-\Delta+V$. Then the corresponding resonance of
the perturbed operator $H+\e V_1$
$$
\lam = \lam_0 - 2\nu\sqrt{-\lam_0}\,\e +O(\e^2),
$$
$\nu$ being given by (\ref{nu}) provided the denominator is  non-zero.
\label{th8}
\end{theorem}

The first order correction to a given resonance is computed in
terms of the solution of the corresponding unperturbed problem
similarly to the well-known self-adjoint \EV \ problem situation.
Naturally, (\ref{nu}) involves $A'$ since the problem is non-linear
with respect to $z$.

The above approach essentially uses the symmetry
of the kernel of $A$ in the sense that (\ref{sa}) is
reduced to the self-adjoint problem, which allows us to compute the correction $\nu$
explicitly. However, a similar argument applies in a generic
situation. Suppose, we have a non-self-adjoint analytic family of operators
$\hat{A}(z,\e)$ whose simple \EVs \ $\lam(z,\e)$ are analytic
functions of the two complex variables (see \cite{Kato}, Ch. VII, \S
4). To find the resonance of the operator $\hat{A}$ we put
$\lam(z,\e)=0$; this equation can be resolved as $z=z(\e)$
provided $\frac{\dpa \lam}{\dpa z} \not= 0$ by the complex analogue
of the implicit function theorem. We therefore conclude that
the resonance $z$ depends analytically on the small parameter
$\e$.

Perturbation formulae for a general analytic operator family are
derived in \cite{Ag}, where the so-called resonance--\EV \ connection
is established in a very abstract context. In section 7 the author
introduces an operator whose \EVs \ coincide with the resonances
of the initial operator; then in section 8 the perturbation
formulae are obtained along the lines of \cite{Kato}. We, on the other hand,
deal with the particular case of the \Schrodinger
operator, computing the first order
term. The proof of theorem~\ref{th8}, being fairly
simple, permits one to switch to the self-adjoint problem,
which is advantageous from the practical point of view when one wants to
find the correction $\nu$.

%%%%%%%%%%%%%%%%%%%%%%%%%%%%%%%
\section{Rapidly decaying potentials}
\label{6}

We describe a technique for finding resonances
which depends on rapid decay of the potential instead of
analyticity assumptions. The method only works in one space
dimension. We treat operators acting in
$L^{2}(\R^+)$ subject to generic boundary conditions at $0$; similar methods apply on the whole line.

\subsection{Existence and uniqueness of the resonance solution}
\label{6.1}

Resonances are conventionally defined to be solutions of the  singular boundary value problem
\begin{eqnarray}
f'' - (V(x) + z^2) f = 0, \qquad x \in\R^+ \ \
\label{eq} \\
a f(0) + b f'(0) = 0, \qquad |a|+|b|>0,
\label{bc} \\
f(x) = \exg + o(|\exd|), \qquad x\to\infty \ \ \ \,
\label{gr}
\end{eqnarray}
Here $z$ is a complex  parameter and
$V(x)$ is a complex-valued function decaying rapidly enough at infinity:
\begin{equation}
\int_{\R^+} |V(x)|\rme^{2z x} \rmd x < \infty.
\label{dec1}
\end{equation}

\begin{lemma}
For all $z$ such that $\Re\,z > 0$ there exists
exactly one solution of (\ref{eq}) satisfying
(\ref{gr}).
\end{lemma}

\begin{proof}
Denote
$$
g(x) =  f(x) - \exg,
$$
then
$$
g'' - (V(x) + z^2) g = V(x)\exg.
$$
Let $u_-(x),\,u_+(x)$ be two linearly independent
solutions of (\ref{eq}) defined by
\begin{equation}
u_\pm(x) \sim d_\pm \rme^{\pm z x}, \qquad x \to \infty,
\label{u12}
\end{equation}
$d_\pm$ being some constants.

We look for the function $g(x)$ in the form
$$
g(x) = C_-(x) u_-(x) + C_+(x) u_+(x),
$$
where the functions $C_-,\,C_+$ solve
$$
C'_- u_- + C'_+ u_+ = 0,
$$
$$
C'_- u'_- + C'_+ u'_+ = V(x)\exg.
$$
In other words,
\begin{equation}
C_\pm(x) =  \mp \frac{1}{W}\int_{x}^{\infty} V(t) \rme^{z t}
u_\mp(t) \rmd t.
\label{c12}
\end{equation}
Here the Wronskian $W = \left| \begin{array}{cc}
u_- & u_+ \\
u'_- & u'_+
\end{array}
\right|$ is constant. Assumptions (\ref{dec1}) and (\ref{u12})
guarantee the convergence of the integrals in (\ref{c12})
and imply
$$
|C_+(x)| = o(\rme^{-2 z x}), \qquad
|C_-(x)| = o(1), \qquad x \to \infty.
$$
This, in turn, means that
\begin{equation}
g(x)=o(\exd), \qquad x \to \infty
\label{z}
\end{equation}
as required.
The way we construct this solution also implies its
uniqueness.
\end{proof}

Now that the lemma is proved we are led
to the following problem: find such values of $z$ that
the equation (\ref{eq}) has a solution satisfying both of
the boundary conditions (\ref{bc}) and (\ref{gr}).
The values of $\lam=-z^2$ are the resonances of our problem.

\subsection{First method}
\label{6.2}

Given any $z$ we can compute at $x=0$
the value of $g(x;z)$ satisfying (\ref{z}).
The resonances are determined by
\begin{equation}
a (1 + g(0;z)) + b (z + g'(0;z)) = 0.
\label{res}
\end{equation}
The resonance solution $f(x;z)$ then satisfies
the boundary condition (\ref{bc}).

One possible method of evaluating $g(0),\,g'(0)$ is as
follows. Denote
$$
\alp_\pm(x) = u'_\pm(x)/u_\pm(x) \mp z.
$$
These two functions
can be easily computed as the solutions of the Cauchy
problems
\begin{eqnarray}
\alp_\pm' + \alp_\pm^2 \pm 2\alp_\pm z - V = 0, \label{ivp} \\
\alp_-(\infty) = 0, \qquad \alp_+(0) = \alp_0.
\label{ic}
\end{eqnarray}
Here the solutions of the problems (\ref{ivp}),~(\ref{ic}) are assumed to be
bounded on $(0,\infty)$. We then have
$$
u_-(x) = u_-(0)\exp\left(\int_{0}^{x} \alp_-(t) \rmd t -
z x\right)
$$
and, hence
\begin{equation}
C_+(x) = - W^{-1}u_-(0)\exp\left(\int_{\R^+} \alp_-(t)
\rmd t\right)I_+(x),
\label{c2}
\end{equation}
where
$$
I_+(x) = \int_{x}^{\infty} V(t) \exp\left(-\int_{t}^{\infty}
\alp_-(\ta) \rmd\ta\right) \rmd t.
$$
Similarly,
\begin{equation}
C_-(x) = W^{-1}u_+(0) I_-(x)
\label{c1}
\end{equation}
where
$$
I_-(x) = \int_{x}^{\infty} V(t) \exp\left(2z t +
\int_{0}^{t} \alp_+(\ta) \rmd\ta\right) \rmd t.
$$

A simple calculation gives us
$$
W = u_-(0)u_+(0)(\alp_+(0) - \alp_-(0) +2z).
$$

We can now evaluate $g(0)$ and $g'(0)$ using (\ref{c2}),
(\ref{c1}). The final formulae are
\begin{equation}
g(0) = \frac{I_-(0) - I_+(0)\exp(\int_{\R^+}
\alp_-(x) \rmd x)}{\alp_0 - \alp_-(0) +2z}\,, \label{z0}
\end{equation}
\begin{equation}
g'(0) = \frac{I_-(0)(\alp_-(0)-z) - I_+(0)(\alp_0+z)
\exp(\int_{\R^+}
\alp_-(x) \rmd x)}{\alp_0 - \alp_-(0) +2z}\,.
\label{z'0}
\end{equation}
We did not specify the initial value $\alp_+(0)$ in (\ref{ivp});
now it is clear that we can take any $\alp_0 \not= \alp_-(0) -
2z$. Computing $g(0),\,g'(0)$ by formulae (\ref{z0}),
(\ref{z'0}) for a range of $z$ we then solve the equation
(\ref{res}) by a standard iterative procedure and, finally
find the resonances.

Following this procedure we have to
integrate (\ref{ivp}), (\ref{ic}) numerically. The equation for
$\alp_-$ is solved from right to left, for $\alp_+$ from
left to right, which means the integration is stable in both cases.
As we compute $\alp_\pm$, at each integration step the
integrals $\int_{0}^{x} \alp_+(t) \rmd t$,
$\int_{x}^{\infty} \alp_-(t) \rmd t$ and
$I_-(x),\,I_+(x)$ are successively evaluated.
Note that the integral $I_-$ may converge rather slowly
depending on the rate of decay of the potential $V(x)$ and on
the size of $|z|$. One has to choose a suitable
method for evaluating $I_-$ according to the behaviour
of the integrand.

\subsection{Second method}
\label{6.3}

Let us normalise the formerly introduced
$u_-(x)$ so that $u_-(x)\exg \to 1$ as $x \to \infty$
($d_-=1$ in (\ref{u12})).
Then
\begin{equation}
u_-(x)=\exd \gam(x), \qquad
\gam(x)=\exp\left(\int_{\infty}^{x} \alp_-(t)\rmd t\right)\,.
\label{u-}
\end{equation}
A general solution of (\ref{eq}) is given by
\begin{equation}
f(x) =  u_-(x)\left(c_1 + c_2 \int_{0}^{x}
\frac{\rmd t}{u_-(t)^2}\right)\,,
\label{gen}
\end{equation}
where $c_1,\,c_2$ are arbitrary constants.
We take $c_2=2z$ and look for such $c_1$ that
$f(x) = \exg + o(u_-(x)),\,x\to\infty$. From (\ref{gen})
we then have
$$
c_1 = \lim_{x\to\infty} \left(\frac{\exg}{u_-(x)} -
2 z \int_{0}^{x} \frac{\rmd t}{u_-(t)^2}\right) =
$$
$$
= u_-(0)^{-1} + \int_{\R^+}
\frac{\exg(z u_-(x) - u'_-(x)) - 2z}{u_-(x)^2} \rmd x =
$$
$$
= u_-(0)^{-1} + \int_{\R^+}
\frac{\rme^{2z x}(\gam(x)(2z -\alp_-(x)) -
2z)}{\gam(x)^2} \rmd x.
$$
Everything is expressed in terms of $\alp_-$ in the
above formula; $u_-(0),\,\gam(x)$ are evaluated by means of (\ref{u-}).
Having found $c_1$ we observe that $f(0) = c_1u_-(0)$,
$f'(0)=c_1u'_-(0) + 2z u_-(0)^{-1}$.
It remains to find the resonances from (\ref{bc}) as before.

The two methods described above along with the method of \cite{AD} based on complex scaling
have been successfully applied to several problems to be considered in
the next section. For a number of problems dealt with here (see the next section)
the numerical efficiency of the two methods appears to be approximately the same.
Choosing between them is a matter of convenience, although there may be situations
when one of them is somewhat better than the other. We observe that the
methods of this section are equivalent within the range of their
consistency.

%%%%%%%%%%%%%%%%%%%%%%%%%%%%%%%
\section{Examples}
\label{7}

\subsection{Gaussian potential}
\label{7.1}

In certain cases no useful information about
the resonances of the operator $H$
can be obtained by the complex scaling method. For example if we put
\begin{equation}
H_gf(x):=-\lap f(x)-\rme^{-x^{2}}f(x)
\label{hg}
\end{equation}
acting in $L^{2}(\R^+)$ subject to Dirichlet  or Neumann boundary
conditions at $0$,  then complex scaling can only be applied for
angles $\theta \leq \pi/2$. However all known
resonances $\lam$ of this operator have $-\pi<\arg(\lam)<-\pi/2$, so
they cannot be computed using the complex scaling
technique \cite{Seidentop}.

We take the operator $H_g$ as an example
and find its resonances  by means of the method of subsection 6.2. The
results given in table~\ref{t1} are obtained up to the accuracy of $10^{-8}$.
Some of the smaller resonances of the same operator, presented very approximately in
\cite{Seidentop}, are at least in a qualitative agreement with our results.
The labels (D) and (N) in table~\ref{t1} relate to Dirichlet and Neumann
boundary conditions respectively.

{\small
\begin{table}
\caption{Resonances for the Gaussian potential}
\label{t1}
\begin{tabular}{|c|c|}
\hline
 $n$  &$\lam_n$ \\ \hline
 1 (D) &-0.52964412\\
 2 (N) &$-1.23692584 -3.48426766 i$\\
 3 (D) &$-1.33981054 -6.70063305 i$\\
 4 (N)& $-1.45338641 -9.95215592 i$\\
 5 (D) &$-1.59459136 -13.14969428 i$\\
 6 (N)   & $-1.69788643 -16.29120850 i$\\
 7 (D)  & $-1.75150063 -19.43849866 i$\\
 8 (N)  & $-1.80493057 -22.61135566 i$\\
 9 (D)  & $-1.87573378 -25.77287400 i$\\
 10 (N)  & $ -1.93288884 -28.91106441 i$\\
 11 (D)  & $-1.96680887 -32.05477086 i$\\
 12 (N)  & $ -2.00314399 -35.21335337 i$\\
 13 (D)  & $ -2.05016204 -38.36530926 i$\\
 14 (N)  & $ -2.13212307 -41.49620539 i$\\
 15 (D)  & $ -2.11478618 -44.64676494 i$\\
 16 (N)  & $ -2.14162017 -47.79467138 i$\\
 17 (D)  & $ -2.16811153 -50.85272557 i$\\ \hline
\end{tabular}
\end{table}
}

\underbar{Note} The operator $H_g$ is also known to have real
\EVs. A simple variational computation shows that there is at
least one \EV \ $\lam_0 \approx -0.335$ corresponding to
the Neumann boundary condition at $x=0$.
On the other hand, it follows from \cite[p.~91]{simon}
that $H_g$ cannot have more than one
negative \EV. This implies that the first entry in table~\ref{t1}
is indeed a resonance with a negligibly small imaginary part.

As well as the Gaussian potential we consider its perturbations
$$
V_\e(x) = -\exp\left(\e^{-1}(1-(1+2\e x^2)^{1/2})\right),
$$
for small $\e >0$.
The reason for making this rather complicated choice of
perturbation is that although one cannot use complex scaling for
the original potential for angles bigger than $\pi/2$, the
approximating potentials allow one to use complex scaling with
angles up to $\pi$, which is enough to capture the resonances of
this operator. The convergence of the first two resonances of the perturbed operator to
those of $H_g$ as $\e\to 0$ is shown in table~\ref{t2}. The
table provides a numerical confirmation of the result of theorem~\ref{th8}.
The resonances of the perturbed and the original operator
agree with the perturbation formulae of section~\ref{5}.
The difference is seen to be of order $O(\e)$ as expected.

{\small
\begin{table}
\caption{Perturbed resonances}
\label{t2}
\begin{tabular}{|c|c|c|}
\hline
$\e$  & $\lam^{(1)}_\e-\lam^{(1)}$& $\lam^{(2)}_\e-\lam^{(2)}$
\\ \hline
$5\cdot10^{-5}$ & $1\cdot10^{-4}$ & $-2\cdot10^{-4} + 6\cdot10^{-4} i$\\
$1\cdot10^{-4}$ & $2\cdot10^{-4}$ & $-3\cdot10^{-4} + 1.1\cdot10^{-3} i$\\
$5\cdot10^{-4}$ & $1\cdot10^{-3}$ & $-1.6\cdot10^{-3} + 5.3\cdot10^{-3} i$\\
$1\cdot10^{-3}$ & $2\cdot10^{-3}$ & $-3.2\cdot10^{-3} + 1.06\cdot10^{-2} i$\\ \hline
\end{tabular}
\end{table}
}

\subsection{A slowly decaying modified Gaussian potential}
\label{7.2}

The operator
\begin{equation}
H_bf(x):=-\lap f(x)+x^2\rme^{-x^{2}/b^2}f(x),
\label{hb}
\end{equation}
$b$ being a real constant, has been studied in \cite{AD}.
Here we provide more detailed information on its resonances for $b=10$.
Corollary 8 is not applicable to the
operator $H_b$ but it seems from the numerical evidence that a
threshold exists. The resonances computed numerically are shown in
figure 1. We used the methods of section 6 to find some of them and
compared the results with those obtained in \cite{AD} by means of complex
scaling. The methods agree with each other for $\Re\,z < 40$; for higher
resonances complex scaling appears to be the best out of the three methods.

Theorem~\ref{6}  can be applied to obtain a bound on all
resonances within the sector
\[
\{ z:-\pi/2 <\arg(z)<0\}.
\]
For this example we have
\[
a(\theta)=b^2\sup\lbrace
s\rme^{-s\cos(\theta)}\sin(2\theta-s\sin\theta):\,
s\in\R^+\rbrace
\]
for every $0\leq\theta <\pi/2$.
This may be used to compute the envelope of the
set $S$ using the parametric formulae (\ref{pf}); this is compared to the actual
resonances (found by complex scaling) in figure 1. Of course $a(\theta)$
may be evaluated purely numerically in more general cases.
Naturally, the bounds on the resonances do not compare closely
with their actual position, but it should be noted that the bounds
are obtained with much less computational effort and do provide
useful qualitative and quantitative information.

\begin{figure}[h]
\begin{picture}(100,230)(0,0)
\includegraphics{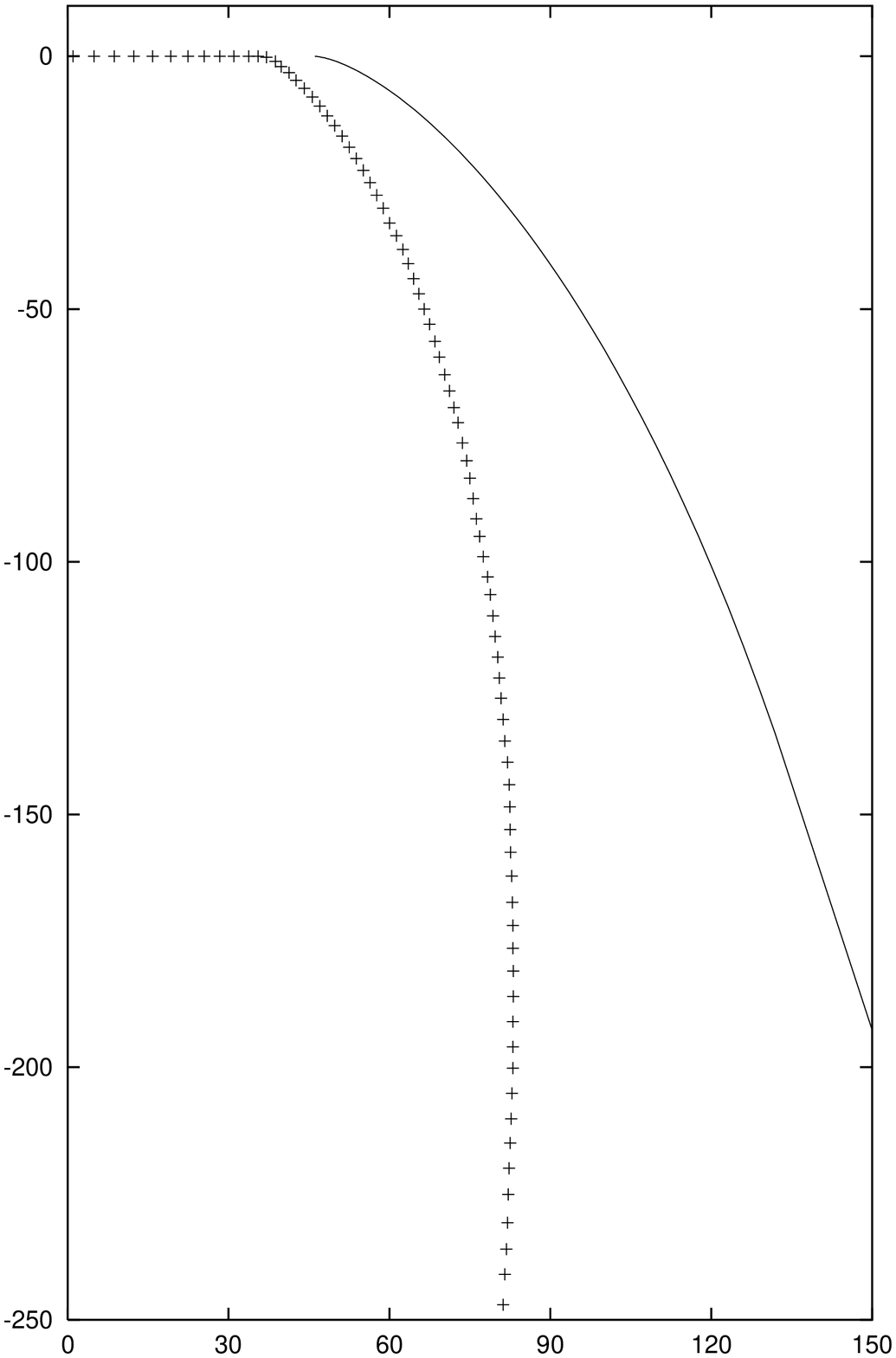}
\end{picture}
\caption{Resonances of $H_b,\,b=10,$ and the set $S$}
\end{figure}

\subsection{The potential of Rittby et al.}
\label{7.3}

We finally study the potential
$$
V(x)=(x^2-J) \rme^{-0.1 x^2} + J, \qquad J=1.6.
$$
In \cite{Br,KLM,REB} one can find controversial remarks on the related
resonance problem and two quite different series of resonances
for the same operator.
Our methods give broadly the same resonances as presented
in \cite{Br}. It appears that the boundary condition at infinity taken in
\cite{KLM} is different from (\ref{gr}) as has been mentioned in
\cite{REB}, hence the difference between the results.

By means of the methods of section~\ref{6} we have computed about 20
resonances, which coincide with those presented in \cite{Br} up to 4 decimal places.
Computing higher resonances by either of the methods of
section~\ref{6} turns out to be quite
difficult, whereas the method of complex scaling still works well.
Indeed, using the technique of \cite{AD} we obtain
stable results which agree very well with all the resonances quoted in
\cite{Br}, whose method is similar.
This suggests that if a potential decays fairly slowly and the complex
scaling method is applicable, then it may often give more accurate
results than the methods of section 6, especially for higher resonances.
We have already commented that for rapidly decaying potentials
like Gaussian, the methods of section~\ref{6} can be efficiently applied.
It is possible to improve those methods
by using an advanced technique for computing the highly oscillatory integrals
involved (see the concluding remark of subsection~\ref{6.2}).

Apart from the resonances found in \cite{Br,REB} we have also
discovered 6 more resonances: two near each of the points
$\lam_1=0.69-7.91i,\,\lam_2=1.26-8.51i,\,\lam_3=2.08-11.61i$.
Each pair includes an even and an odd resonance very close to one another.
This closeness hinders accurate computations and is probably the reason why
those resonances are missing in \cite{Br,REB}. None of our three methods
allows us to locate them accurately, but there is clear
evidence that they exist.

\bibliographystyle{plain}
\bibliography{art,book}

\end{document}